\documentclass[11pt]{article}
\usepackage[dvips]{graphicx}

\usepackage{amsfonts,amsmath, amssymb, latexsym, boxedminipage}
\def\squarebox#1{\hbox to #1{\hfill\vbox to #1{\vfill}}}
\def\qed{\hspace*{\fill}
        \vbox{\hrule\hbox{\vrule\squarebox{.667em}\vrule}\hrule}\smallskip}
\newenvironment{proof}{\begin{trivlist}
  \item[\hspace{\labelsep}{\em\noindent Proof.~}]
  }{\qed\end{trivlist}}
\newtheorem{lemma}{Lemma}[section]
\newtheorem{theorem}[lemma]{Theorem}
\newtheorem{corollary}[lemma]{Corollary}
\newtheorem{proposition}[lemma]{Proposition}

\newtheorem{observation}[lemma]{Observation}
\newtheorem{definition}[lemma]{Definition}

\newcommand{\comment}[1]{}
\def\squareforqed{\hbox{\rlap{$\sqcap$}$\sqcup$}}
\def\qed{\ifmmode\squareforqed\else{\unskip\nobreak\hfil
\penalty50\hskip1em\null\nobreak\hfil\squareforqed
\parfillskip=0pt\finalhyphendemerits=0\endgraf}\fi}
\newcommand{\power}{\mbox{$\mathbb P$}}

\newcommand{\field}{\mbox{$\mathbb F$}}
\newcommand{\nats}{\mbox{$\mathbb N$}}
\newcommand{\ints}{\mbox{$\mathbb Z$}}

\newcommand{\ind}{\mbox{\rm ind}}
\newcommand{\Gdag}{\vec{G}}

\newcommand{\maxset}{\max\{\Gdag\}}

\newlength{\tablength}
\newlength{\spacelength}
\newcommand{\tabstar}{\hspace*{\tablength}}
\newcommand{\spacestar}{\hspace*{\spacelength}}
\def\obeytabs{\catcode`\^^I=\active}
{\obeytabs\global\let^^I=\tabstar}
{\obeyspaces\global\let =\spacestar}
\newenvironment{display}{\begingroup\obeylines\obeyspaces\obeytabs}{\endgroup}
\newenvironment{prog}{\begin{display}\parskip0pt\sf}{\end{display}}

\title{Vertex coloring acyclic digraphs and their corresponding hypergraphs\thanks{Some of the results presented here appeared in a preliminary form in \cite{agtive}.}}

\author{Geir Agnarsson
\thanks{Department of Mathematical Sciences,
George Mason University,
MS 3F2,
4400 University Drive,
Fairfax, VA 22030,
{\tt geir@math.gmu.edu}.}
\and \'{A}g\'{u}st S.~Egilsson
\thanks{Science Institute,
University of Iceland 
IS-107 Reykjavik, Iceland,
{\tt egilsson@hi.is}.}
\and Magn\'{u}s M.~Halld\'{o}rsson
\thanks{Department of Computer Science, 
University of Iceland,
IS-107 Reykjavik, Iceland,
{\tt mmh@hi.is}.}
}

\date{}

\begin{document}

\maketitle

\begin{abstract}
\noindent We consider vertex coloring of an acyclic digraph $\Gdag$ in
such a way that two vertices which have a common ancestor in $\Gdag$
receive distinct colors. Such colorings arise in a natural way
when bounding space for various genetic
data for efficient analysis. We discuss the corresponding
{\em down-chromatic number} and derive an upper bound as a function 
of $D(\Gdag)$, the maximum number of descendants of a given vertex, and the 
degeneracy of the corresponding hypergraph. Finally we determine an 
asymptotically tight upper bound of the down-chromatic number
in terms of the number of vertices of $\Gdag$ and $D(\Gdag)$.

\vspace{3 mm}
%% 05C15 - coloring, 05C20 - digraphs, 05C85 - algorithms,
%% 05C90 - applications, 05C65 - hypergraphs, 06A05 - posets


\vspace{2 mm}

\end{abstract}

\section{Introduction}
\label{sec:intro}

%% \subsection{Purpose}

%% Subject of this paper
The purpose of this article is to discuss a special kind 
of vertex coloring for acyclic digraphs, where vertices 
with a common ancestor must receive distinct colors. 
We discuss some properties of such colorings, 
similarity and differences with strong hypergraph colorings,
and derive an upper bound which, in addition, yields an
efficient coloring procedure. 

Digraphs representing various biological phenomena and knowledge
are ubiquitous in the life sciences and in drug discovery
research, e.g.~the gene ontology digraph maintained by the Gene
Ontology Consortium~\cite{GO}. An overview of several projects relating 
to indexing of semistructured data (i.e.~acyclic digraphs) can be found 
in~\cite{Ab-Bu-Su}. In these biological digraphs it is important
to be able to access the ancestors of nodes in a fast and
efficient manner.

% Problem considered: Fast DB access to ancestor set of a node
Consider the problem of finding a representation an acyclic
digraph in a database to allow for fast access to the set of 
ancestors of a given node. The ancestors of a node are its 
in-neighbors in the transitive closure of the digraph.
If the digraph is sparse and shallow, the transitive closure is also
sparse. Thus, an adjacency matrix representation would be neither efficient
nor fast. On the other hand, a matrix has the advantage
of corresponding nicely to the relational representation of modern databases.
In a database relation that corresponds to an adjacency matrix, 
the non-empty elements in each column correspond to the in-neighbors
of the node indexing that column. This vertex set can
then be combined by joins with other tables that are also indexed by vertices,
giving an effective language of querying based on graphical properties.
Thus, it would be preferable to find a representation that both has a
matrix structure, yet consists of relatively small rows.

% A compacted matrix repres. of digraphs
One compact matrix representation would be to store the adjacency lists 
in compacted array form, where a list with $k$ elements is stored in the
first $k$ array elements. 
In this case, however, there is no easy way of accessing all the edges
entering a given vertex. While this could be alleviated by storing the
inverted adjacency matrix, note that in the context of database
access, rows are conceptually different from columns.
Instead, we seek a compacted representation
where all in-edges of a given node are stored in the same column.
We say that a many-to-one mapping of vertices to columns that preserves
adjacency lists, has the {\em AC-property}.
By recording the mapping of nodes to their respective column storing
their in-neighbors, one obtains the same desirable properties of
adjacency matrices in the context of a relational database. If the
graph is sparse, the possibilities of storage reduction are significant.
The actual improvement is related to the number of colors needed in a
certain coloring of the digraph, which we now briefly discuss.

% Def. of down-coloring + Relationship with compact repres. of TC of digraph
A proper \emph{down-coloring} of a digraph is a vertex coloring where
vertices with a common ancestor receive different colors. The
\emph{down-chromatic number} of a digraph is the minimum number of
colors in a down-coloring of the digraph. In a compacted matrix
representation of the transitive closure of a digraph, we assign multiple
vertices to the same column, but in such a way that their in-adjacency
lists must be disjoint. Two vertices have disjoint sets of in-neighbors in
the transitive closure if, and only if, they have no common ancestor.
Therefore, a down-coloring of a digraph corresponds to a valid compacted
representation of its transitive closure, and the down-chromatic number
is the minimum number of columns needed in such a representation.

\vspace{3 mm}

\noindent {\sc Example:} Consider the digraph $\Gdag$, on $n=6$ vertices
representing genes, where a directed edge from one vertex to 
a second one indicates that the first gene is an ancestor of the
second gene.
\begin{eqnarray*}
V(\Gdag) & = & \{g_1,g_2,g_3,g_4,g_5,g_6\}, \\
E(\Gdag) & = & \{ (g_1,g_4),(g_1,g_5),(g_2,g_4),(g_2,g_6),(g_3,g_5),(g_3,g_6)\}.
\end{eqnarray*}
In the adjacency matrix representation of this 6 node digraph, we assign a column 
to each vertex $g_i$. As we see in the left diagram of Table~\ref{tab:n=6}, 
most of the entries of this $6\times 6$ matrix are empty. 
\begin{table}
\begin{center}
\begin{tabular}{c|cccccc}
          & $g_1$ & $g_2$ & $g_3$ & $g_4$ & $g_5$ & $g_6$  \\
\hline
$g_1$     & 1     & 0     & 0     & 1     & 1     & 0      \\
$g_2$     & 0     & 1     & 0     & 1     & 0     & 1      \\
$g_3$     & 0     & 0     & 1     & 0     & 1     & 1      \\
$g_4$     & 0     & 0     & 0     & 1     & 0     & 0      \\
$g_5$     & 0     & 0     & 0     & 0     & 1     & 0      \\
$g_6$     & 0     & 0     & 0     & 0     & 0     & 1      \\
\end{tabular}
\ \ \ $\rightarrow$ \ \ \ 
\begin{tabular}{c|cccccc}
          & 1     &  2    &  3    \\
\hline
$g_1$     & $g_1$ & $g_4$ & $g_5$ \\
$g_2$     & $g_6$ & $g_4$ & $g_2$ \\
$g_3$     & $g_6$ & $g_3$ & $g_5$ \\
$g_4$     & --    & $g_4$ & --    \\
$g_5$     & --    & --    & $g_5$ \\
$g_6$     & $g_6$ & --    & --    \\
\end{tabular}
\end{center}
\caption{An example for $n=6$.}
\label{tab:n=6}
\end{table}
By reducing the number of columns in such a way
that the AC-property still holds, we obtain a smaller and more compact 
$6\times 3$ matrix representation as seen on the right diagram of Table~\ref{tab:n=6}.
There, the $i$-th row still contains the descendants of $g_i$ and 
the ancestors of $g_i$ are those $g_j$'s whose rows $g_i$ appears in.
Note that, (i) each $g_i$ appears in exactly one column and, (ii) two
genes appear in the same column only if their sets of ancestors are
disjoint. Here we can view the column numbers 1,\ 2, and 3 as distinct colors
assigned to each vertex. Note further that in this example the transitive closure
of $\Gdag$ is simply $\Gdag$ itself.
An explicit example of how such a coloring can speed up queries 
in the gene ontology digraph can be found in the appendix of~\cite{Egilsson}.

Hence, the following two questions regarding such colorings, one 
computational and the other theoretical, are quite natural:
(1) For a given digraph (of no particular structure!) how can 
we assign reasonably few colors to the vertices/columns
efficiently, and (2) in general, how large can the discrepancy theoretically
be between the actual minimum number of colors needed and
the obvious lower bound of needed colors?

\subsection*{Our Results}
The contributions of this paper are threefold. First, we establish a
close link between down coloring digraphs and strong coloring hypergraphs.
Second, we give efficiently computable bounds on the down chromatic
number in terms of the inductiveness of the related hypergraph and 
$D(\Gdag)$, the maximum number of descendants of a given vertex.
And thirdly, we give a tight bound on the discrepancy between the down
chromatic number and the lower bound $D(\Gdag)$. This also has
independent interest as characterizing the largest ratio of the
strong chromatic number of a hypergraph to the sum of the number of
edges and the number of vertices.

\subsection*{Related Work}
Some special classes of such acyclic digraphs are studied 
in~\cite{GeirAgust}, in particular those of height two
in which every vertex has an in-degree of two. 
For a brief introduction and additional references to the ones mention here, 
we refer to~\cite{GeirAgust}. 

Note that acyclic digraphs are often called 
directed acyclic graphs or DAG's by computer
scientists, as is the case in~\cite[p.~194]{Shaffer}.  

A straight forward condition of a vertex coloring of a digraph $\Gdag$
is to insist that two vertices $u$ and $v$ receive distinct color
if there is a directed edge from $u$ to $v$ in $\Gdag$.
Such a coloring is, of course, the same as coloring the vertices of the underlying 
graph $G$ of $\Gdag$ (by forgetting the orientation of the directed edges) 
in the usual sense. 

Another vertex coloring of digraphs that relies on the direction of
the edges is the {\em dichromatic number} of a digraph $\Gdag$,
as studied in~\cite{Jacob} and~\cite{Su}, which is defined as 
the minimum number of colors needed to vertex color $\Gdag$
in such a way that no monochromatic directed cycle is created.

Strong colorings of hypergraphs have been
studied, but not quite to the extent of various other types of 
colorings of hypergraphs. Since strong colorings are generalizations 
of the usual vertex colorings of graphs, the determination of the exact
strong chromatic number of a hypergraph is in general a 
daunting task. Most results in this direction in the literature on
strong colorings are restricted to some very special types of hypergraphs.
In~\cite{Weifan} a nice survey of various aspects of hypergraph 
coloring theory is found, containing almost all 
fundamental results in the past three decades.
What we are concerned here is not necessarily
an {\em exact} computation of the strong chromatic
number, but rather a good theoretical upper bound
that is valid for all possible corresponding digraphs
$\Gdag$. In Section~\ref{sec:discrep} however, we discuss the
asymptotics of how large the exact down-chromatic number can be.

\section{Basic definitions}
\label{sec:defs}

We attempt to be consistent with standard graph theory notation
in~\cite{DWest}, and the notation in~\cite{Trotter} when applicable. 
%% Digraphs
For a natural number $n\in\nats$ we let $[n] = \{1,\ldots,n\}$.
A {\em simple digraph} is a finite simple directed graph
$\Gdag = (V,E)$, where $V=V(\Gdag)$ is a finite set
of vertices and $E=E(\Gdag)\subseteq V\times V$ is a set of
directed edges. The digraph $\Gdag$ is said
to be {\em acyclic} if $\Gdag$ has no directed cycles.
Henceforth $\Gdag$ will denote an acyclic digraph
in this section.
%% Posets
The binary relation $\leq$ on $V(\Gdag)$ defined by
\begin{equation}
\label{eqn:poset}
u \leq v \Leftrightarrow u = v, \mbox{ or there is a directed path
from $v$ to $u$ in $\Gdag$,}
\end{equation}
is reflexive, antisymmetric and transitive
and therefore a partial order on $V(\Gdag)$.
Hence, whenever we talk about $\Gdag$ as a {\em poset}, the partial
order will be the one defined by (\ref{eqn:poset}). 
The {\em transitive closure} of $\Gdag$ is the poset $\Gdag$ 
viewed as a digraph, that is the digraph $\Gdag^*$ on $V(\Gdag)$ where 
$(v,u)\in E(\Gdag^*)$ iff $u < v$. 
By the {\em height} of $\Gdag$ as a poset, we mean the number
of vertices in the longest directed path in $\Gdag$.
We denote by $\max\{\Gdag\}$ the set of maximal vertices of $\Gdag$ 
with respect to the partial order $\leq$.
%% Down-sets
For vertices $u,v\in V(\Gdag)$ with 
$u\leq v$, we say that $u$ is a {\em descendant} of $v$, and
$v$ is an {\em ancestor} of $u$. 
The {\em closed principal down-set} or simply the {\em down-set} $D[u]$ of a vertex $u\in V(\Gdag)$ 
is the set of descendants of $u$
in $\Gdag$, that is, $D[u] = \{x\in V(\Gdag) : x \leq u\}$.
Likewise, the {\em open principal down-set} or the {\em open down-set} of 
a vertex $u$ is $D(u) = D[u]\setminus \{u\}$.
\begin{definition}  
\label{def:down-coloring}
A {\em down-coloring} of $\Gdag$ is a map
$c:V(\Gdag) \rightarrow [k]$ satisfying
\[
u,v \in D[w] \mbox{ for some $w\in V(\Gdag)$ }
\Rightarrow c(u)\neq c(v)
\] 
for every $u,v\in V(\Gdag)$.
The {\em down-chromatic number} of $\Gdag$,
denoted by $\chi_d(\Gdag)$, is the least $k$
for which $\Gdag$ has a proper down-coloring
$c:V(\Gdag) \rightarrow [k]$.
\end{definition}
Clearly, in an undirected graph $G$ the vertices in a clique
must all receive distinct colors in a proper vertex coloring
of $G$. Therefore $\omega(G)\leq \chi(G)\leq |V(G)|$ where 
$\omega(G)$ denotes the clique number of $G$.
Similarly, if $D(\Gdag) = \max_{u\in V(\Gdag)}\{|D[u]|\}$
for our acyclic digraph $\Gdag$, we clearly have 
$D(\Gdag)\leq \chi_d(\Gdag) \leq |V(\Gdag)|$.
Hence, when considering down-colorings, it can be useful to map the problem to
one on undirected graphs. Given an acyclic digraph $\Gdag$, 
the corresponding simple undirected {\em down-graph} $G'$ has the same
set of vertices, with each pair
of vertices connected that are contained in the same principal down-set:
\begin{eqnarray*}
V(G') & = & V(\Gdag), \\
E(G') & = & \{ \{u,v\} : 
u,v\in D[w]\mbox{ for some }w\in V(\Gdag)\}.
\end{eqnarray*}
In this way we have transformed the problem of down-coloring
the digraph $\Gdag$ to the problem of vertex coloring
the simple undirected graph $G'$ in the usual sense,
and we have $\chi_d(\Gdag) = \chi(G')$.
Hence, from the point of down-colorings, both
$\Gdag$ and $G'$ are equivalent.

As observed in~\cite[Obs.~2.3]{GeirAgust} we have: 
\begin{observation}
\label{obs:no-func}
There is no function $f : \nats\rightarrow\nats$ with
$\chi_d(\Gdag) \leq f(D(\Gdag))$ for all acyclic digraphs $\Gdag$.
\end{observation}
However, although not a function of $D(\Gdag)$ alone, there
are computable parameters such that $\chi_d(\Gdag)$ can
be bounded by functions in terms of these parameters. That will
be the purpose of the following section.

\section{Hypergraph representations}
\label{sec:reps}

In this section we discuss alternative representations
of our digraph $\Gdag$, and define some parameters which
we will use to bound the down-chromatic number $\chi_d(\Gdag)$.

We first consider the issue of the height of digraphs.
We say that two digraphs on the same set of vertices
are {\em equivalent} if every down-coloring of
one is also a valid down-coloring of the other, that is, if they induce
the same undirected down-graph.
We show that for any acyclic digraph 
$\Gdag$ there is an equivalent acyclic digraph
$\Gdag_2$ of height two with $\chi_d(\Gdag) = \chi_d(\Gdag_2)$.
\begin{lemma}
\label{lmm:height2}
Any down-graph $G'$ of an acyclic digraph $\Gdag$ is also
a down-graph of an acyclic digraph $\Gdag_2$ of height two.
\end{lemma}
\begin{proof}
The derived digraph $\Gdag_2$ has the same vertex set as $\Gdag$, 
while the edges all go from $\maxset$ to $V(\Gdag)\setminus\maxset$,
where $(u,v) \in E(\Gdag_2)$ if, and only if, $v \in D(u)$.
In this way we see that two vertices in $\Gdag$ have a 
common ancestor if, and only if, they have a common ancestor
in $\Gdag_2$. Hence, we have the proposition.
\end{proof}
Therefore, when considering down-colorings of digraphs, we can 
by Lemma~\ref{lmm:height2}
assume them to be of height two. 

Recall that a {\em hypergraph} is $H$ is set system
on $V$, that is $H = (V, \mathcal{E})$ where $V$ is a set of vertices
and $\mathcal{E}$ is a set (possibly a multiset) of 
subsets of $V$ called {\em hyperedges}. A hypergraph is
{\em simple} if $\mathcal{E}$ is not a proper multiset (that is,
$\mathcal{E}\subseteq \power(V)$, the power set of $V$), 
and each hyperedge
has cardinality 2 or more. For a given hypergraph $H$, simple or not,
denote by $V(H)$ the set of its vertices and $\mathcal{E}(H)$ the
set of its hyperedges. Two vertices of a hypergraph $H$ are 
{\em neighbors} in $H$ if they are contained in the
same edge in $\mathcal{E}(H)$.  An edge in $\mathcal{E}(H)$ containing 
just one element is called {\em trivial}. The largest cardinality of a 
hyperedge of $H$ will be denoted by $\sigma(H)$. 
To every simple hypergraph $H$ there is
an associated simple {\em clique graph} $G$ on the same vertices as $H$
where two vertices are connected iff they are contained in the same
hyperedge. Note that two distinct simple hypergraphs can have  
identical clique graphs.

There is a natural correspondence between acyclic digraphs and 
certain hypergraphs.
\begin{definition}
\label{def:down-hyp}
For a digraph $\Gdag$, the corresponding {\em down-hypergraph} $H_{\Gdag}$
of $\Gdag$ is defined by:
\begin{eqnarray*}
 V(H_{\Gdag})           & = & V(\Gdag)\setminus\maxset, \\
 \mathcal{E}(H_{\Gdag}) & = & \{ D(u) : u \in \maxset\}.
\end{eqnarray*}
Conversely, for a hypergraph $H$ the corresponding {\em up-digraph} $\Gdag_H$ 
of $H$ is defined by:
\begin{eqnarray*}
          V(\Gdag_H)   & = & V(H)\cup \{w_e : e\in \mathcal{E}(H)\}, \\
          E(\Gdag_H)   & = & \{ (w_e,u) : u\in e\in \mathcal{E}(H)\}.
\end{eqnarray*}
\end{definition}
Note that with the notation from above we have for any digraph $\Gdag$
with no isolated vertices
that $\Gdag_{H_{\Gdag}} = \Gdag_2$, the equivalent digraph of height two
from here above. We summarize in the following:
\begin{observation}
\label{obs:uep}
For any hypergraph $H$ we have $H_{\Gdag_H} = H$ and for any
digraph $\Gdag$ of height two with no isolated vertices we
have $\Gdag_{H_{\Gdag}} = \Gdag$.
\end{observation}
Hence, for our down-coloring purposes, digraphs are equivalent
to digraphs of height two with no isolated vertices, which then again are equivalent to
hypergraphs, where vertices in the same hyperedge receive
different colors. This is precisely a {\em strong coloring} of a hypergraph $H$, that 
is a map $\Psi : V(H) \rightarrow [k]$ such that 
$u,v\in e$ for some $e\in \mathcal{E}(H)$, implies $\Psi(u)\neq\Psi(v)$.
The {\em strong chromatic number} $\chi_s(H)$ is the least number 
$k$ of colors for which $H$ has a proper strong coloring
$\Psi : V(H) \rightarrow [k]$.
Just as for graphs, when considering strong colorings of hypergraphs,
we can, with no loss of generality, restrict to simple hypergraphs.

For an acyclic digraph $\Gdag$ we see that an optimal strong coloring
of $H_{\Gdag}$ will yield and optimal down-coloring of $\Gdag$, simply
by completing the colorings of $\maxset$ in a greedy fashion. In the 
case where $\chi_s(H_{\Gdag}) = \sigma(H_{\Gdag})$, then since 
$D(\Gdag) = \sigma(H_{\Gdag}) + 1$, we have 
$\chi_d(\Gdag) = \chi_s(H_{\Gdag}) + 1$. Otherwise, when
$\chi_s(H_{\Gdag}) > \sigma(H_{\Gdag})$, we always have at least
one available color from the set $\{1,2,\ldots, \chi_s(H_{\Gdag})\}$
to complete the down-coloring of $\Gdag$ in a legitimate and optimal 
fashion. Hence we have $\chi_d(\Gdag) = \chi_s(H_{\Gdag})$ in this case.
We summarize in the following.
\begin{theorem}
\label{thm:G-HG}
For an acyclic digraph $\Gdag$ we have
\[
\chi_d(\Gdag) = \left\{ \begin{array}{ll}
    \chi_s(H_{\Gdag}) + 1 & \mbox{ if } \chi_s(H_{\Gdag}) = \sigma(H_{\Gdag}), \\
    \chi_s(H_{\Gdag})     & \mbox{ if } \chi_s(H_{\Gdag}) > \sigma(H_{\Gdag}). \\ 
                        \end{array}
                \right.
\]
\end{theorem}
We can also characterize the down chromatic number precisely in terms
of the strong chromatic number of related hypergraph. The \emph{closed
  down hypergraph} $\hat{H}_{\Gdag}$ has the same vertex set as
$H_{\Gdag}$ but the edgeset $E(\hat{H}_{\Gdag}) = \{D[u] : u \in \maxset\}$.

\begin{observation}
For an acyclic digraph $\Gdag$, we have $\chi_d(\Gdag) = \chi_s(\hat{H}_{\Gdag})$.
\end{observation}
The down-graph $G'$ of $\Gdag$ is precisely the clique-graph of the
closed down-hypergraph $\hat{H}_{\Gdag}$.

\subsection*{Computable bounds}

For a hypergraph $H = (V(H),\mathcal{E}(H))$ the {\em degree} $d_H(u)$, 
or just $d(u)$, of a vertex $u\in V(H)$ is the number of non-trivial edges containing $u$.
The minimum and maximum degree of $H$ are given by 
$\delta(H) = \min_{u\in V(H)}\{d_H(u)\}$ and 
$\Delta(H) = \max_{u\in V(H)}\{d_H(u)\}$ respectively. 
The subhypergraph $H[S]$ of $H$, induced by a set $S$ of vertices, is 
given by 
\begin{eqnarray*}
V(H[S]) & = & S, \\
\mathcal{E}(H[S]) & = & \{ X\cap S : X \in \mathcal{E}(H)\mbox{ and }|X \cap S| \ge 2 \}.
\end{eqnarray*}
\begin{definition}
\label{def:ind-hyp}
Let $H$ be a simple hypergraph. The {\em degeneracy} or the {\em inductiveness}  
of $H$, denoted by $\ind(H)$, is given by
\[
\ind(H) = \max_{S\subseteq V(H)}\left\{ \delta(H[S])\right\}.
\]
If $k\geq \ind(H)$, then we say that $H$ is {\em $k$-degenerate}
or {\em $k$-inductive}.
\end{definition}
Note that Definition~\ref{def:ind-hyp} is a generalization
of the degeneracy or the inductiveness of a usual undirected graph $G$,
given by $\ind(G) = \max_{H\subseteq G}\left\{ \delta(H)\right\}$. 
Note that the degeneracy of a (hyper)graph 
is always greater than or equal to the degeneracy of
any of its sub(hyper)graphs.

To illustrate, let us for a brief moment discuss the degeneracy 
of an important class of simple graphs, namely that of
simple planar graphs. Every subgraph of a simple planar graph 
is again planar. Since every planar graph has
a vertex of degree five or less, the degeneracy
of every planar graph is at most five. This is 
the best possible for planar graphs, since 
the graph of the icosahedron is planar and
5-regular. That a planar graph has degeneracy
of five, implies that it can be vertex colored
in a simple greedy fashion with at most six
colors. The degeneracy has also been used to
bound the chromatic number of the square $G^2$
of a planar graph $G$, where $G^2$ is a graph obtained
from $G$ by connecting two vertices of $G$ if, and only
if, they are connected in $G$ or they have a common 
neighbor in $G$ (see \cite{GeirMagnus}.)
In general, the degeneracy of an undirected graph
$G$ yields an ordering
$\{ u_1,u_2,\ldots ,u_n\}$ of $V(G)$, such that 
each vertex $u_i$ has at most $\ind(G)$ neighbors
among the previously listed vertices 
$u_1,\ldots, u_{i-1}$. Such an ordering provides
a way to vertex color $G$ with 
at most $\ind(G) + 1$ colors in an efficient
greedy way, and hence we have in general
that $\chi(G)\leq \ind(G) + 1$. 

The degeneracy of a simple hypergraph is also
connected to a greedy vertex coloring of it,
but not in such a direct manner as for a regular undirected
graph, since, as noted, the number of neighbors of a given
vertex in a hypergraph is generally much larger than its
degree. 
\begin{theorem}
\label{thm:ind-hyp}
If the simple undirected graph $G$ is the clique
graph of the simple hypergraph $H$ then 
$\ind(G) \leq \ind(H)(\sigma(H) - 1)$.
\end{theorem}
\begin{proof}
For each $S\subseteq V(G)=V(H)$, let $G[S]$ and $H[S]$ be
the subgraph of $G$ and the subhypergraph of $H$ induced
by $S$, respectively. Note that for each $u\in S$, each 
hyperedge in $H[S]$ which contains $u$, has at most 
$\sigma(H[S])-1\leq\sigma(H)-1$ other vertices 
in addition to $u$. By definition of $d_{H[S]}(u)$, we
therefore have that 
$d_{G[S]}(u) \leq d_{H[S]}(u)(\sigma(H)-1)$, and hence
\begin{equation}
\label{eqn:d-u}
\delta(G[S]) \leq \delta(H[S])(\sigma(H)-1).
\end{equation}
Taking the maximum of (\ref{eqn:d-u}) among 
all $S\subseteq V(G)$ yields the theorem.
\end{proof}
Recall that the {\em intersection graph} of a collection $\{A_1,\ldots,A_n\}$
of sets, is the simple graph with vertices $\{u_1,\ldots,u_n\}$,
where we connect $u_i$ and $u_j$ if, and only if, $A_i\cap A_j\neq \emptyset$.

Directly by definition of the inductiveness we have the
following.
\begin{observation}
\label{obs:ind-tree}
For a simple connected hypergraph $H$, then $\ind(H)=1$ if, and only if,
the intersection graph of its hyperedges $\mathcal{E}(H)$ is a tree.
\end{observation}
What Observation~\ref{obs:ind-tree} implies, is that
edges of $H$ can be ordered as $\mathcal{E}(H) = \{e_1,\ldots,e_m\}$,
such that each $e_i$ intersects exactly one edge from
the set $\{e_1,\ldots, e_{i-1}\}$. If now $G$ is the clique graph
of $H$, this implies that $\ind(G)=\sigma(H)-1$ and hence
$\chi(G) = \sigma(H)$. Therefore, by Theorem~\ref{thm:G-HG}, we have 
in general the following:
if $\ind(H_{\Gdag}) = 1$, then $\chi_s(H_{\Gdag}) = \chi(G) = \sigma(H_{\Gdag})$,
and hence $\chi_d(\Gdag) = \sigma(H_{\Gdag}) + 1 = D(\Gdag)$.
Otherwise, if $\ind(H_{\Gdag}) > 1$, then by Theorem~\ref{thm:ind-hyp} 
we have
\[
\chi_s(H_{\Gdag}) = \chi(G) \leq \ind(G) + 1\leq \ind(H_{\Gdag})(\sigma(H_{\Gdag})-1) + 1.
\]
Since now $D(\Gdag) = \sigma(H_{\Gdag}) + 1$ we have therefore
the following corollary.
\begin{corollary}
\label{cor:summary}
If $\Gdag$ is an acyclic digraph, then its down-chromatic number 
satisfies the following:
\begin{enumerate}
  \item If $\ind(H_{\Gdag}) = 1$ then $\chi_d(\Gdag) = D(\Gdag)$.
  \item If $\ind(H_{\Gdag}) > 1$ then $\chi_d(\Gdag) \leq \ind(H_{\Gdag})(D(\Gdag) - 2) + 1$.
\end{enumerate}
Moreover, in both cases the given upper bound of colors 
can be used to down-color $\Gdag$ in an efficient greedy
fashion.
\end{corollary}
{\sc Example:}
Let $k,m\in\nats$, let $A_1,\ldots, A_k$ be disjoint sets,
each $A_i$ containing exactly $m$ vertices. Let $H(k,m)$
be the hypergraph with 
\begin{eqnarray*}
V(H(k,m))           & = & \bigcup_{i\in [k]}A_i, \\
\mathcal{E}(H(k,m)) & = & \{A_i\cup A_j : i\neq j, \ \ \{i,j\}\subseteq [k]\},
\end{eqnarray*}
Let $\Gdag(k,m) = \Gdag_{H(k,m)}$ be the up-digraph of the hypergraph $H(k,m)$.
Clearly $\Gdag(k,m)$ is a simple acyclic digraph
on $km + {k\choose 2}$ vertices and with
${k\choose 2}\cdot 2m = k(k-1)m$ directed edges.
Further, $\sigma(H(k,m)) = 2m$ and so $D(\Gdag(k,m)) =  2m + 1$.
Since each vertex is contained in exactly $k$ hyperedges we have
$\ind(H(k,m))=k-1$. Hence, by Corollary~\ref{cor:summary}, we
obtain that $\chi_d(\Gdag(k,m))\leq (k-1)(2m-1)+1 = \Theta(km)$,
which agrees with the asymptotic value
of the actual down-chromatic number $km$ (also a $\Theta(km)$ function).
Hence, up to a constant (of 2), Corollary~\ref{cor:summary} is asymptotically tight.

\section{Discrepancy between parameters}
\label{sec:discrep}

So far we have discussed how to approximate the down-chromatic
number $\chi_d(\Gdag)$ of an acyclic digraph $\Gdag$ in terms of $D(\Gdag)$ and
$\ind(H_{\Gdag})$, the inductiveness of the corresponding down-hypergraph.
In this section we will discuss the relative discrepancy between $D(\Gdag)$
and the actual down chromatic number $\chi_d(\Gdag)$, and determine
a tight asymptotic upper bound for their ratio.

If $H(k,m)$ is the hypergraph defined in the last example of the previous section, then
for $\Gdag_{H(k,m)}$ we clearly have 
\[
\frac{\chi_d( \Gdag_{H(k,m)})}{D(\Gdag_{H(k,m)})} = \frac{km}{2m+1}\rightarrow \infty
\]
as $k\rightarrow \infty$ and $m$ is fixed. Hence, allowing an unbounded
number of vertices of $\Gdag$, the above ratio clearly can become arbitrarily large
even when $D(\Gdag(k,m)) =  2m + 1$ is fixed.

The purpose of this last section is to derive a tight upper bound for 
$\chi_d( \Gdag)/D(\Gdag)$
among all acyclic digraphs $\Gdag$ with $D(\Gdag)$ bounded and with
bounded number of vertices. 
\begin{definition}
\label{def:rel-dis}
For $n,\delta \in \nats$ define the {\em relative down-coloring discrepancy}, or
simply the {\em rdcd}, $d_d(\delta,n)$
by 
\[
d_d(\delta,n) = \max_{|V(\Gdag)| \leq n,\ D(\Gdag)\leq \delta}
\left\{
\frac{\chi_d(\Gdag)}{D(\Gdag)}
\right\},
\]
where the maximum is among all acyclic digraphs $\Gdag$ satisfying the 
stated conditions.
\end{definition}
Note that for a hypergraph $H$, it holds that
$D(\Gdag_H) = \sigma(H) + 1$ and $|V(\Gdag_H)| = |V(H)| + |\mathcal{E}(H)|$.
Hence, for $n,\sigma \in \nats$ we define the {\em relative strong-coloring discrepancy}, or
simply the {\em rscd}, $d_s(\sigma,n)$
by 
\begin{equation}
\label{eqn:hyp}
d_s(\sigma,n) = \max_{|V(H)| + |\mathcal{E}(H)|\leq n,\ \sigma(H) \leq \sigma}
\left\{
\frac{\chi_s(H)}{\sigma(H) + 1}
\right\},
\end{equation}
where the maximum is taken among all hypergraphs $H$.
By Lemma~\ref{lmm:height2} and Observation~\ref{obs:uep} we 
have the following.
\begin{observation}
\label{obs:disc}
For $n,\delta, \sigma \in \nats$ we have $d_d(\sigma + 1,n) = d_s(\sigma,n)$.
\end{observation}
Although our original motivation for the relative discrepancy $d_d(\delta,n)$ is
given by Definition~\ref{def:rel-dis}, by Observation~\ref{obs:disc}
it suffices to (and in some ways is more natural to) 
determine a tight upper bound of $d_s(\sigma,n)$ for given $n,\sigma \in \nats$
from (\ref{eqn:hyp}).
\begin{definition}
\label{def:r+}
For $n,\sigma \in \nats$, let $r^+(\sigma,n)$ denote the positive root
of the quadratic polynomial $x + \frac{x(x-1)}{\sigma(\sigma-1)}=n$
in terms of $x$.
\end{definition}
Using Definition~\ref{def:r+} we now can state our first theorem.
\begin{theorem}
\label{thm:r+}
For $n,\sigma \in \nats$ the rscd $d_s(\sigma,n)$ satisfies
\[
d_s(\sigma,n) \leq \frac{r^+(\sigma,n)}{\sigma + 1}
\]
\end{theorem}
\begin{proof}
Let $H$ be a hypergraph with $|V(H)| + |\mathcal{E}(H)|\leq n$,
$\chi_s(H) = x\in \nats$ and $\sigma(H) = \sigma$.
In this case there is an optimal strong $x$-coloring of the vertices
of $H$. Let $V_1,\ldots,V_x$ be corresponding partition 
of $V(H)$ into color classes. For each $i$ and $j$ with $1\leq i<j\leq x$,
there is at least one hyperedge $e_{i\/j}\in \mathcal{E}(H)$ that
contains one vertex from $V_i$ and one vertex from $V_j$. Since
$|e|\leq \sigma$ for each $e\in \mathcal{E}(H)$, each hyperedge
can cover at most $\binom{\sigma}{2}$ sets of two vertices that are
colored by distinct pairs of colors. Since there are $\binom{x}{2}$
pairs of colors, the number of hyperedges of $H$ must satisfy
\[
|\mathcal{E}(H)| \geq \frac{\binom{x}{2}}{\binom{\sigma}{2}} = \frac{x(x-1)}{\sigma(\sigma-1)}.
\]
Since each color class $V_i$ is nonempty, we must have 
$|V(H)|\geq \chi_s(H) = x$. Combining the last two inequalities
we obtain, in particular, that
\begin{equation}
\label{eqn:hyp-cond}
n \geq |V(H)| + |\mathcal{E}(H)| \geq x + \frac{x(x-1)}{\sigma(\sigma-1)}.
\end{equation}
Viewing $n$ and $\sigma$ as arbitrary but fixed, we obtain by (\ref{eqn:hyp}) 
and (\ref{eqn:hyp-cond}) that 
\begin{equation}
\label{eqn:slick}
d_s(\sigma,n)\leq 
\max_{x + \frac{x(x-1)}{\sigma(\sigma-1)} \leq n}\left\{\frac{x}{\sigma + 1}\right\}.
\end{equation}
By solving the corresponding
quadratic inequality $x + \frac{x(x-1)}{\sigma(\sigma-1)} \leq n$ in terms
of $x$, keeping in mind that $x$ is positive, we have that
$0<x<r^+(\sigma,n)$. The maximum value of the fraction 
$x/(\sigma + 1)$ is clearly taken when $x$ is at maximum, 
that is for $x = r^+(\sigma,n)$. Hence, by (\ref{eqn:slick}) we obtain
\[
d_s(\sigma,n) \leq \frac{r^+(\sigma,n)}{\sigma + 1},
\]
which completes the proof.
\end{proof}
By solving the quadratic equation $x + x(x-1)/(\sigma(\sigma-1)) = n$ for $x$ we get
\[
r^+(\sigma,n) = \frac{\sqrt{\sigma(\sigma-1)n}}
{\sqrt{1 + \frac{(\sigma(\sigma-1)-1)^2}{4\sigma(\sigma-1)n}} 
+ \frac{\sigma(\sigma-1)-1}{\sqrt{4\sigma(\sigma-1)n}}},
\]
so it is immediate that $r^+(\sigma,n)\leq \sqrt{\sigma(\sigma-1)n}$ and 
\begin{equation}
\label{eqn:r+assymp}
\lim_{n\rightarrow\infty}\frac{r^+(\sigma,n)}{\sqrt{n}} = \sqrt{\sigma(\sigma-1)}.
\end{equation}
Hence, by Theorem~\ref{thm:r+} we obtain the following.
\begin{corollary}
\label{cor:rel-dis}
For $n,\sigma \in \nats$ the rscd $d_s(\sigma,n)$ satisfies
\[
d_s(\sigma,n) \leq \frac{\sqrt{\sigma(\sigma-1)}}{\sigma + 1}\sqrt{n}.
\]
\end{corollary}
Note that for a fixed $\sigma$ the upper bounds for $d_s(\sigma,n)$ 
in Theorem~\ref{thm:r+} and Corollary~\ref{cor:rel-dis} are by (\ref{eqn:r+assymp})
asymptotically the same as $n\rightarrow \infty$.
We now argue that the upper bound from Corollary~\ref{cor:rel-dis} is 
asymptotically tight in the sense that
\begin{equation}
\label{eqn:limsup}
\limsup_{ n\rightarrow\infty}\frac{d_s(\sigma,n)}{\sqrt{n}} = \frac{\sqrt{\sigma(\sigma-1)}}{\sigma + 1}
\end{equation}
for infinitely many values of $\sigma$. More specifically, we will
show that there is an infinite collection $(\sigma_i)_{i\geq 1}$ such that
for each $i$ there is again an infinite collection
$(n_{i\/j})_{j\geq 1}$ with the property that there exists a hypergraph $H_{i\/j}$
with $|V(H_{i\/j})| + |\mathcal{E}(H_{i\/j})| = n_{i\/j}$ and $\sigma(H_{i\/j}) = \sigma_i$ 
that matches the upper bound of Theorem~\ref{thm:r+}, that is
\[
\frac{\chi_s(H_{i\/j})}{\sigma(H_{i\/j}) + 1} = \frac{r^+(\sigma_i,n_{i\/j})}{\sigma_i + 1}
\]
for each $i$ and $j$. This together with Theorem~\ref{thm:r+} and Corollary~\ref{cor:rel-dis}
will yield (\ref{eqn:limsup}). For this we need some additional terminology for hypergraphs.

Recall that a {\em balanced incomplete block design} or a {\em BIBD} for short,
is a simple hypergraph $H = (V,\mathcal{B})$ where $V$ is
a finite set of vertices and $\mathcal{B}\subseteq \power(V)$ 
is a collection of hyperedges where 
(i) all the hyperedges have
the same cardinality that is strictly less than that of $V$,
(ii) each vertex is contained in the same $r>0$ number of hyperedges,
and  
(iii) each pair of vertices is contained in exactly $\lambda>0$ 
hyperedges. In this case the vertices are sometimes called
{\em varieties} and the hyperedges {\em blocks}.
In a series of three papers~\cite{WilsonI}, \cite{WilsonII} and 
\cite{WilsonIII}, R.~M.~Wilson proved 
that for any given $k,\lambda\in \nats$ there
exists a constant $C = C(k,\lambda)\in \nats$ such that
for any $v\geq C$ satisfying (i) $\lambda(v-1)\equiv 0\bmod(k-1)$
and (ii) $\lambda v(v-1)\equiv 0 \bmod k(k-1)$, then there
exists a BIBD on $v$ vertices, where each hyperedge has cardinality
$k$ and where each vertex is contained in $\lambda$ hyperedges.
In particular, for $\lambda = 1$, we have with our notation and 
terminology from above the following.
\begin{corollary}
\label{cor:Wilson}
For each $\sigma\in \nats$ there exists a constant $K = K(\sigma)\in \nats$ 
such that for all $k\geq K$ with $k^2\equiv k\bmod\sigma$, 
we have
\[
d_s\left(\sigma,k(\sigma - 1 + k) + 1 - \frac{k(k-1)}{\sigma}\right) = 
\frac{k(\sigma-1) + 1}{\sigma + 1}.
%\frac{\left\lfloor r^+\left(\sigma,k(\sigma - 1 + k) + 1 - \frac{k(k-1)}{\sigma}
%\right)\right\rfloor}{\sigma + 1},
\]
%where $r^+(\sigma,n)$ is as in Definition~\ref{def:r+}.
\end{corollary}
\begin{proof}
By R.~M.~Wilson, there is a $C = C(\sigma)\in \nats$ such that for all
$v\geq C$ satisfying $v-1\equiv 0\bmod(\sigma-1)$ and 
$v(v-1)\equiv 0\bmod\sigma(\sigma-1)$, there is a BIBD, call it $H$, on $v$
vertices such that each hyperedge has exactly $\sigma$ vertices,
each pair of vertices is contained in exactly one hyperedge.
In particular, the number of hyperedges of $H$ is given by
\[
|\mathcal{E}(H)| =  \frac{\binom{v}{2}}{\binom{\sigma}{2}} = \frac{v(v-1)}{\sigma(\sigma-1)}.
\]
The conditions on $v$ mean that $v = k(\sigma-1)+1$ where $k^2\equiv k\bmod\sigma$
and $k\in\nats$ is large enough, say $k \geq K = K(\sigma)$.
Since each pair of the $v$ vertices is contained in a hyperedge, we clearly
have $\chi_s(H) = v$. We also clearly have $\sigma(H) =  \sigma$ and hence, 
in this case we have 
$\chi_s(H)/(\sigma(H) + 1) = v/(\sigma + 1) =( k(\sigma - 1) + 1)/(\sigma+1)$. Also,
for $v = k(\sigma-1)+1$ we have
\[
n = v + \frac{v(v-1)}{\sigma(\sigma-1)} = k(\sigma-1+k)+1 -\frac{k^2-k}{\sigma}.
\]
By (\ref{eqn:hyp}) this implies that
\[
d_s\left(\sigma,k(\sigma - 1 + k) + 1 - \frac{k(k-1)}{\sigma}\right)
\geq \frac{k(\sigma-1) + 1}{\sigma + 1}.
\]
Since $r^+\left(\sigma,k(\sigma - 1 + k) + 1 - \frac{k(k-1)}{\sigma}\right) =  k(\sigma-1)+1$,
we have by Theorem~\ref{thm:r+} the corollary.
\end{proof}
We conclude this section by an explicit and self contained 
construction of a class of hypergraphs for which the asymptotic value in 
(\ref{eqn:limsup}) can also be reached. First note that
for $\sigma\in\nats$ and $v = \sigma^k$ then
$v-1\equiv 0\bmod(\sigma-1)$ and $v(v-1)\equiv 0\bmod\sigma(\sigma-1)$
hold for all $k\in\nats$.

Before proving Proposition~\ref{prp:affine} we recall some notations
and results: Every finite field has cardinality of a prime power $q = p^n$.
If $\ints_p$ denotes the integers modulo $p$, then the unique field 
$\field_q$ of cardinality $q$ can be given as the quotient 
$\field_q = \ints_p[X]/(X^q-X)$ which turns out to be the splitting field 
of the polynomial $X^q-X$ over $\ints_p$. In particular, $\field_r$ is a 
subfield of $\field_q$ whenever $r= p^m$
and $m\leq n$ (See~\cite[p.~278]{Hung}.) 
The {\em affine $d$-space} over a field $F$ is a tuple 
$(F^d, \mathcal{L})$ where $F^d$ consists of all ordered $d$-tuples
$\tilde{x} = (x_1,\ldots,x_d)$ where each $x_i\in F$, and where $\mathcal{L}$
is the collection of all lines $\{ \tilde{a} + t\tilde{b} : t\in F\}\subseteq F^d$.
In particular, the {\em affine plane}
over a field $F$ is the affine $2$-space over $F$, which is a BIBD 
with $\lambda = 1$ (See~\cite[p.~199]{LintWilson}.)
\begin{proposition}
\label{prp:affine}
If $\sigma = p^k$ is a prime power, then for any $m\in\nats$ there
is a BIBD on $\sigma^m$ vertices, where each hyperedge has cardinality
$\sigma$ and where each pair of vertices is contained in exactly one hyperedge.
\end{proposition}
\begin{proof}
For a given prime $p$ and positive integers $k$ and $m$, let 
$H = (V(H),\mathcal{E}(H))$ be the affine $m$-space over 
the field $\field_{\sigma}$ where $\sigma = p^k$. Then $H$
is a BIBD on $|V(H)| = {\sigma}^m$ vertices where each hyperedge
(i.e.~line) has cardinality $\sigma = p^k$ and where each pair
of vertices is contained in exactly one hyperedge. This completes the proof.
\end{proof}
{\sc Remarks:} (i) The condition that each pair of vertices is contained
in exactly one hyperedge is a natural geometric condition called the 
{\em Euclid's first postulate}, when vertices are viewed as points and
hyperedges as lines. (ii) Note that the affine plane over any field $F$,
in particular for the finite field $\field_q$, satisfies
the {\em Euclidean parallel postulate} (aka {\em Euclid's fifth postulate}), 
that for any three vertices,
not all contained in a hyperedge, there is precisely one
hyperedge containing the third vertex, that is disjoint from
the unique hyperedge containing the first two vertices.
However, for $d\geq 3$ the affine $d$-space
has the {\em hyperbolic parallel property}: for any three vertices,
not all contained in the same hyperedge, there are two or
more hyperedges containing the third vertex that are also
disjoint from the unique hyperedge containing the first two
vertices. (See~\cite{Greenberg}.) (iii) It is a well-known conjecture, 
whether or not there exists an affine plane on $n$ vertices
when $n$ is not a power of a prime, is still open. (See~\cite{Bogart}.)

From Proposition~\ref{prp:affine} we deduce the following
corollary.
\begin{corollary}
\label{cor:affine}
If $\sigma = p^k$ is a prime power, then for any $m\in\nats$ we have
\[
d_s\left(\sigma, \sigma^m + \frac{\sigma^{m-1}(\sigma^m - 1)}{\sigma-1}\right)
=\frac{\sigma^{m}}{\sigma + 1}.
\]
\end{corollary}
\begin{proof}
We clearly have
\[
r^+\left(\sigma, \sigma^m + \frac{\sigma^{m-1}(\sigma^m - 1)}{\sigma-1}\right) = \sigma^m,
\]
and hence by Theorem~\ref{thm:r+} we have
\[
d_s\left(\sigma, \sigma^m + \frac{\sigma^{m-1}(\sigma^m - 1)}{\sigma-1}\right)
\leq \frac{r^+(\sigma,n)}{\sigma + 1} = \frac{\sigma^m}{\sigma + 1},
\]
yielding the upper bound.

By Proposition~\ref{prp:affine} there is a BIBD $H$ on $\sigma^m$ vertices, 
where each hyperedge has cardinality $\sigma$ and where each pair of vertices is 
contained in exactly one hyperedge. In this case we have
$\chi_s(H) = \sigma^m$ and $\sigma(H) = \sigma$ and hence 
\[
d_s\left(\sigma, \sigma^m + \frac{\sigma^{m-1}(\sigma^m - 1)}{\sigma-1}\right)
\geq \frac{\chi_s(H)}{\sigma(H) + 1} = \frac{\sigma^m}{\sigma + 1},
\]
yielding the lower bound and so the proof is complete.
\end{proof}
{\sc Remark:} Throughout this article we have assumed our 
digraphs to be acyclic. However, we note 
that the definition of down-coloring can be easily extended
to a regular cyclic digraph $\Gdag$ by interpreting the notion of
descendants of a vertex $u$ to mean the set of nodes reachable from $u$.
In fact, if $\Gdag$ is an arbitrary digraph, then there is an 
equivalent acyclic digraph $\Gdag'$, on the same set of vertices,
with an identical down-graph:
First form the {\em condensation}
$\hat{G}$ of $\Gdag$ by shrinking each strongly connected 
component of $\Gdag$ to a single vertex.
Then form $\Gdag'$ by replacing each node of $\hat{G}$ 
which represents a strongly connected component of $\Gdag$ on a set 
$X\subseteq V(\Gdag)$ of vertices, with an arbitrary vertex $u\in X$,
and then add a directed edge from $u$ to each 
$v\in X \setminus \{u\}$. This completes the construction. 

Observe that each node $v\in X$ has exactly the same neighbors in the
down-graph of $\Gdag'$ as $u$, as it is a descendant of $u$ and $u$
alone.  Further, if node $v$ was in a different strong component of
$\Gdag$ than $u$ but was reachable from $u$, then it will continue to
be a descendant of $u$
in $\Gdag'$. Hence, the down-graphs of $\Gdag$ and $\Gdag'$ are identical.

\subsection*{Acknowledgments}

The authors are grateful to John L.~Pfaltz for his interest and
encouragements to write this article. Also, sincere thanks to 
Jim Lawrence for pointing out the very relevant work of 
Richard M.~Wilson on BIBDs.

\today
\end{document}